\documentclass[10pt]{article}
\usepackage[cp1251]{inputenc}
\usepackage[russian]{babel}
\begin{document}

\centerline{\bf Задача о книжной полке}

\centerline{Л.Радзивиловский, Г. Юргин
\footnote{Исследование поддержано грантом РФФИ 14-01-00548.}}

\bigskip
В задачнике "Математического просвещения" номер 18 была

\medskip
{\small
{\bf Задача 8.} На полке в некотором порядке стоят тома, пронумерованные числами от 1 до $n$. Библиотекарь берёт том, стоящий не на своём месте, и ставит его на правильное место; при этом некоторые тома сдвигаются.

(а) Докажите, что процесс перестановки томов остановится. ({\it Фольклор})

(б) Постарайтесь получить оценку на число тактов этого процесса, например полиномиальную. ({\it А.Я.Белов})
}

\medskip

Пункт (а) имеет много решений. Каждое решение, явно или косвенно, даёт оценку на количество тактов. Доказательство того, что оценка
точная, может быть построено на процессе, состоящем из большого количества ходов.
Мы докажем, что количество ходов экспоненциально зависит от $n$, а именно, что процесс заканчивается
не более чем за $2^n-2^{\frac{n}{2}}$
тактов, но может продолжаться $2^{n-1}-1$ тактов.
Многие наши знакомые пытались доказать полиномиальную оценку, но это, как теперь выяснилось,  невозможно. Между $2^n-2^{\frac{n}{2}}$  и $2^{n-1}-1$  остаётся зазор; мы приглашаем читателей к дальнейшим размышлениям.

\medskip
Начнём с {\bf примера}. Для двух книг, стоящих в неправильном порядке, процесс закончится за 1 ход. Можно
считать, что мы при этом не снимаем книгу 2 с полки.

Продолжим строить примеры по индукции. Предположим, что при начальном расположении книг в порядке $n,1,2,\dots,n-1$  можно упорядочить книги, сделав ровно $2^{n-1}-1$  ходов, и при этом не нужно снимать том $n$ с полки.

Возьмём $n+1$
книг,
расположенных в порядке $n+1,1,2,\dots,n$. Сделаем ряд ходов, не снимая с полки тома 1 и $n+1$. При этом том  $n+1$ будет оставаться на месте, и вся игра будет происходить на $n$ правых позициях. Можно считать, что всего есть $n$  томов, причём том 2 играет
роль тома 1, том 3 играет роль тома 2, \dots, том $n$  играет роль тома $n-1$, а том 1 играет роль тома $n$
$$n+1,2,3,4,\dots,n,1.$$
Все тома стоят на своих местах, кроме первого и последнего. Теперь первый том переставляется на своё место, и получается
$$1,n+1,2,3,4,\dots,n.$$
Теперь имеется
$n$  книг, с уже знакомой нам начальной ситуацией, и можно произвести ещё   $2^{n-1}-1$ ходов; при этом нет необходимости снимать том $n+1$  с полки.
Всего мы сделали $\left(2^{n-1}-1\right)+1+\left(2^{n-1}-1\right)=2^n-1$  ходов.

\medskip
Теперь {\bf опишем полуинвариант}, который позволит нам
оценить сверху
количество ходов, а заодно и доказать конечность процесса.

Представим себе две строчки, каждая из $n-1$ лампочек.
В первой строчке первая (слева) лампочка горит тогда и только тогда, когда первый том на месте; вторая лампочка горит тогда и только тогда, когда второй том находится на одном из двух первых мест; 
...;
$(n-1)$-я
лампочка горит тогда и только тогда, когда  $(n-1)$-й том находится на одном из $n-1$ первых мест.

Вторая строчка лампочек устроена по похожему принципу, только наоборот: первая (слева) лампочка горит тогда и только тогда, когда последний том на своём месте; вторая лампочка горит тогда и только тогда, когда предпоследний том на одном из двух последних мест; \dots; ($n-1)$-я лампочка горит тогда и только тогда, когда  $(n-1)$-й с конца том на одном из  $n-1$ последних мест.

Предположим, что какая-нибудь лампочка гаснет. Мы утверждаем, что одновременно с этим зажигается какая-то лампочка левее в той же строчке. Действительно, если $k$-й том вытеснен с первых $k$  мест, это значит, что был поставлен на своё место том $j$, где $j<k$. Аналогично, если $k$-й с конца том был вытеснен с последних $k$  мест, то только потому, что на место был поставлен
$j$-й ($j<k$)
с конца том.

Итак, если считать, что лампочки в каждой строчке представляют собой двоичные цифры какого-то числа, то оба числа не убывают. Более того, с каждым ходом одно из чисел возрастает, ибо постановка тома на своё место всегда зажигает какую-то лампочку (если том сдвинулся влево, то лампочку в первой строке, а если вправо, то во второй).

Таким образом,
у нас есть два двоичных $(n-1)$-значных числа, и с каждым ходом одно из них возрастает. Их сумма не меньше чем $0$  и не больше чем $2^n-2$, поэтому всего можно сделать не более чем $2^n-2$  ходов.

Несложно понять, что эту оценку можно улучшить.
Каждой книге с номером $k,\ 1<k<n$ соответствуют две лампочки --- $k$-я в первой строчке и $(n-k)$-я во второй, и в любой момент горит хотя бы одна из них. Отсюда получается, что начальная сумма двух $(n-1)$-значных чисел не меньше $2^{\frac{n}{2}}-2$, а
 число ходов не превышает  $2^n-2^{\frac{n}{2}}$.

{\small
{\bf Замечание.}
Имеется задача на схожий сюжет.

{\it На полке стоят тома. Разрешается брать  любые  два  тома,
идущие в обратном  порядке,  и  менять  их  местами  (в  обратном
порядке значит, что том с большим  номером  стоит  левее  тома  с
меньшим номером). Доказать, что этот процесс закончится.}
\medskip

Доказательство основано на том, что общее число  {\it беспорядков} (т.е. пар томов, когда том с бОльшим номером стоит левее тома с меньшим номером) уменьшается, отсюда получается полиномиальная  верхняя оценка на число тактов $n(n-1)/2$. Несложно привести оптимальный пример -- тома идут в порядке убывания, переставляются соседи. С другой стороны, за $n-1$ операцию всегда можно осуществить сортировку. Всегда ли можно это сделать за м\'еньшее число операций и на сколько?
}

Авторы признательны А.~Я.~Белову за постановку задачи и внимание к работе.
\bigskip

Л.~Радзивиловский, Тель-Авивский университет

levr78@hotmail.com

Г.~Юргин, Москва, лицей  Вторая школа, 10 класс

g.y.98@mail.ru

\end{document}